\documentclass[11pt,twoside]{amsart}
\usepackage{amscd,amssymb}

\title[local constancy of characters]{On the local constancy of characters}
\author{Jonathan Korman}
\email{jkorman@math.toronto.edu}
\address{The University of Toronto \\ Toronto, Ontario M5S 3G3}
\date{\today}

\begin{document}

%\linespread{1.6}
\theoremstyle{plain}
\newtheorem{thm}{Theorem}
\newtheorem*{thm*}{Theorem}
\newtheorem{lem}[thm]{Lemma}
\newtheorem{cor}[thm]{Corollary}
\newtheorem{propn}[thm]{Proposition}
\newtheorem*{propn*}{Proposition}
\newtheorem{claim}[thm]{Claim}
\newtheorem{eqn}[thm]{eqnarray}

\theoremstyle{definition}
\newtheorem{rmk}[thm]{Remark}
\newtheorem{rmks}[thm]{Remarks}
\newtheorem{defn}[thm]{Definition}
\newtheorem*{defn*}{Definition}
\newtheorem{ex}[thm]{Example}
\newtheorem{hypo}[thm]{Hypothesis}

\begin{abstract}
The character of an irreducible admissible representation of a
$p$-adic reductive group is known to be a constant function in
some neighborhood of any regular semisimple element $\gamma$ in
the group. Under certain mild restrictions on $\gamma$, we give an
explicit description of a neighborhood of $\gamma$ on which the
character is constant.
\end{abstract}

%\begin{flushright} {\LARGE \textmd{Preliminary Version}} \end{flushright}
%\vspace{1cm}

\maketitle

\setcounter{section}{-1}

%\tableofcontents

\section*{Introduction}
Let $k$ be a $p$-adic field of characteristic zero, and let
$\mathbf{G}$ be a connected reductive algebraic group defined over
$k$. We denote by $G$ the group of $k$-rational points of
$\mathbf{G}$, and by $\mathfrak{g}$ the Lie algebra of $G$. Let
$\pi$ be an irreducible admissible representation of $G$, and let
$\Theta_\pi$ be the (distribution) character of $\pi$.
In~\cite{Queens} Harish-Chandra showed that $\Theta_\pi$ can be
represented by a function (also denoted by $\Theta_\pi$) which is
locally integrable on $G$ and locally constant on the set
$G^{reg}$ of regular semisimple elements in $G$. Thus for any
$\gamma \in G^{reg}$ there exists {\it some} neighborhood of
$\gamma$ on which the character is constant. In~\cite[Theorem~2,
p.~483]{Howe}, R. Howe gave an elementary proof of
Harish-Chandra's result for general linear groups. In this paper
we give a precise version of local constancy (near compact regular
semisimple tame elements) for all reductive groups. The outline of
the approach given here follows the elementary argument of Howe.

Let $\mathfrak{g}_{x,r}$ (resp.\ $G_{x,|r|}$) be the Moy-Prasad
lattices~\cite{Moy-Prasad-1} in $\mathfrak{g}$ (resp.\ open
compact subgroups of $G$), normalized as
in~\cite[\S1.2]{Fiona-Kim}. Let $G_{cpt}$ denote the set of
compact elements in $G$.
%(or in the notation of DeBacker $G_0=\cup_{x\in\mathcal{B}(G)} G_{x,0}$)
For a maximal $k$-torus $T$, let $T_r$ denote its filtration
subgroups (section~\ref{root}). Let $\rho (\pi)$ denote the depth
of $\pi$~\cite[\S 5]{Moy-Prasad-1}.

Fix a regular semisimple element $\gamma$ and let
$\mathbf{T}:=C_\mathbf{G}(\gamma)^\circ$ be the connected
component of its centralizer; $\mathbf{T}$ is a maximal $k$-torus
in $\mathbf{G}$. We assume that it splits over some tamely
ramified finite Galois extension $E$ of $k$. Let $T$ denote the
group of $k$-rational points of $\mathbf{T}$. When $\gamma\in T
\cap G_{cpt}$ we attach to it the nonnegative rational number
$s(\gamma)$. Using the filtration subgroups $T_r$ and the
parameter $s(\gamma)$, we characterize a neighborhood of $\gamma$
on which the character $\Theta_\pi$ is constant. Whether or not
this neighborhood of constancy is maximal is not addressed here.

The main result of this paper is the following
(Theorem~\ref{mainresult}).

\begin{thm*}
Let $r = \mathrm{\max} \{ s(\gamma),\rho (\pi)\}+s(\gamma)$. The
character $\Theta_\pi$ is constant on the set ${}^G(\gamma
T_{r+})$.
\end{thm*}

We now give a brief sketch of the proof. Let $K$ be any open
compact subgroup of $G$. Decompose $\Theta_\pi$ into a countable
sum of `partial trace' operators $\Theta_d$, according to the
irreducible representations $d$ of $K$ (see section~\ref{partial
traces}). For $G=GL_n$, Howe proved~\cite[p.~499]{Howe} the
following key fact. If $X$ is a compact subset of $G^{reg}$, then
$\Theta_d$ vanishes on $X$ for all $d$ not in a certain finite set
$F$ (which depends only on $X$). It follows (see proof of
Theorem~\ref{mainresult}), that $\Theta_\pi (f) = \int_X (\sum_{d
\in F} \Theta_d)(x) f(x) dx$ for all $f\in C_c^\infty (X)$. Hence
$\Theta_\pi$ is represented on $X$ by the locally constant
function $\sum_{d \in F} \Theta_d$.

The main part of this paper is concerned with formulating an
analogue of the above key fact for reductive groups (see
Corollary~\ref{cord}).

The rational number $s(\gamma)$, defined in
section~\ref{regular_depth}, is used (Corollary~\ref{cord}) to
make a precise choice of a set $X$ and a subgroup $K$.
Corollary~\ref{cord} characterizes a finite set $F$ of
representations, such that for all $d$ not in $F$, $\Theta_d$
vanishes on $X$ (see Remark~\ref{supp_of d} for the significance
of this fact). Thus the representations $d\in F$ are those which
play a role in understanding the character $\Theta_\pi$ near
$\gamma$. The proof of this corollary relies on a special case
(Corollary~\ref{cor to propn}), in which we only consider
$1$-dimensional $d$. Such representations have an explicit
description in terms of cosets in the lie algebra $\mathfrak{g}$.
In section~\ref{tech_lemma} we develop the technical tools, using
Moy-Prasad lattices, to handle these cosets. Once we have a
characterization of the set $F$, we can make precise statements
about the neighbourhood of constancy of the character near
$\gamma$ (Theorem~\ref{mainresult}).

\subsection*{Acknowledgments.}
I would like to thank Fiona Murnaghan, Jeff Adler, Stephen
DeBacker, Ju-Lee Kim and Joe Repka for helpful conversations and
comments. I would also like to thank the referees for their
careful reading and detailed comments.

\section*{Notation and Conventions}
Let $k$ be a $p$-adic field (a finite extension of some
$\mathbb{Q}_p$) with residue field $\mathbb{F}_{p^{n}}$. Let $\nu$
be a valuation on $k$ normalized such that
$\nu(k^{\times})=\mathbb{Z}$.

For any algebraic extension field $E$ of $k$, $\nu$ extends
uniquely to a valuation (also denoted $\nu$) of $E$.

We denote the ring of integers in $E$ by $R_E$ (write $R$ for
$R_k$), and the prime ideal in $R_E$ by $\wp_E$ (write $\wp$ for
$\wp_k$).

%Fix a uniformizing element $\varpi_E$ (resp.\ $\varpi$) in $\wp_E$ (resp.\ $\wp$).

Let $\mathbf{G}$ be a connected reductive group defined over $k$,
and $\mathbf{G}(E)$ the group of $E$-rational points of
$\mathbf{G}$. We denote by $G$ the group of $k$-rational points of
$\mathbf{G}$. Denote the Lie algebras of $\mathbf{G}$ and
$\mathbf{G}(E)$ by {\boldmath $\mathfrak{g}$} and {\boldmath
$\mathfrak{g}$}$(E)$, respectively. Write $\mathfrak{g}$ for the
Lie algebra of $k$-rational points of {\boldmath $\mathfrak{g}$}.

Let $\mathcal{N}$ be the set of nilpotent elements in
$\mathfrak{g}$. There are different notions of nilpotency, but
since we assume that char$(k)=0$, these notions all coincide.

Let $\mbox{Ad}$ (resp.\ $\mbox{ad}$) denote the adjoint
representation of $\mathbf{G}$ (resp.\ {\boldmath $\mathfrak{g}$})
on its Lie algebra {\boldmath $\mathfrak{g}$}. For elements $g\in
G$ and $X\in\mathfrak{g}$ (resp.\ $x\in G$) we will sometimes
write ${{}^gX}$ (resp.\ ${{}^gx}$) instead of $\mbox{Ad}(g)X$
(resp.\ $gxg^{-1}$). For a subset $S$ of $\mathfrak{g}$ (resp.\
$G$) let ${}^G S$ denote the set $\{^gs\,|\,g\in G \mbox{ and }
s\in S\}$.

Let $n$ denote the (absolute) rank of $\mathbf{G}$. We say that an
element $g\in G$ is {\it regular semisimple} if the coefficient of
$t^n$ in det$(t-1+$Ad$(g))$ is nonzero. We denote the set of
regular semisimple elements in $G$ by $G^{reg}$. Similarly we say
that an element $X\in\mathfrak{g}$ is {\it regular semisimple} if
the coefficient of $t^n$ in det$(t-$ad$(X))$ is nonzero. We denote
the set of regular semisimple elements in $\mathfrak{g}$ by
$\mathfrak{g}^{reg}$. Let $G_{cpt}$ denotes the set of compact
elements in $G$. For a subset $S$ of $G$ we will sometimes write
$S_{cpt}$ for $S\cap G_{cpt}$.

For a subset $S$ of $\mathfrak{g}$ (resp.\ $G$) let $[S]$ denote
the characteristic function of $S$ on $\mathfrak{g}$ (resp.\ $G$).

For any compact group $K$, let $K^\wedge$ denote the set of
equivalence classes of irreducible, continuous representations of
$K$.

Let $\pi$ be an irreducible admissible representation of $G$. We
denote by $\Theta_\pi$ the character of $\pi$ thought of as a
locally constant function on the set $G^{reg}$. Let $\rho(\pi)$
denote the depth of $\pi$~\cite[\S 5]{Moy-Prasad-1}.

\section{Preliminaries}
\subsection{Apartments and buildings.} For a finite extension $E$
of $k$, let $\mathcal{B}(\mathbf{G},E)$ denote the extended
Bruhat-Tits building of $\mathbf{G}$ over $E$; write
$\mathcal{B}(G)$ for $\mathcal{B}(\mathbf{G},k)$. It is known
(e.g.~\cite{Prasad}) that if $E$ is a tamely ramified finite
Galois extension of $k$ then $\mathcal{B}(\mathbf{G},k)$ can be
embedded into $\mathcal{B}(\mathbf{G},E)$ and its image is equal
to the set of Galois fixed points in $\mathcal{B}(\mathbf{G},E)$.
If $\mathbf{T}$ is a maximal $k$-torus in $\mathbf{G}$ that splits
over $E$, let $\mathcal{A}(\mathbf{T},E)$ be the corresponding
apartment over $E$. Let $\mathbf{X}^*(\mathbf{T},E)$ (resp.\
$\mathbf{X}_*(\mathbf{T},E)$) denote the group of $E$-rational
characters (resp.\ cocharacters) of $\mathbf{T}$.

It is known in the tame case~\cite[\S 1.9]{Adler} that there is a
Galois equivariant embedding of $\mathcal{B}(\mathbf{T},E)$ into
$\mathcal{B}(\mathbf{G},E)$, which in turn induces an embedding of
$\mathcal{B}(\mathbf{T},k)$ into $\mathcal{B}(\mathbf{G},k)$. Such
embeddings are only unique modulo translations by elements of
$\mathbf{X}_*(\mathbf{T},k)\otimes \mathbb{R}$, however their
images are all the same and are equal to the set
$\mathcal{A}(\mathbf{T},E)\cap\mathcal{B}(\mathbf{G},k)$. From now
on we fix a $T$-equivariant embedding $i:\mathcal{B}(\mathbf{T},k)
\longrightarrow \mathcal{B}(\mathbf{G},k)$, and use it to regard
$\mathcal{B}(\mathbf{T},k)$ as a subset of
$\mathcal{B}(\mathbf{G},k)$; write $x$ for $i(x)$.

{\bf Notation.} We write $\mathcal{A}(\mathbf{T},k)$ for
$\mathcal{A}(\mathbf{T},E)\cap\mathcal{B}(\mathbf{G},k)$. This is
well defined independent of the choice of $E$~\cite{Yu}. Moreover,
$\mathcal{A}(\mathbf{T},k)$ is the set of Galois fixed points in
$\mathcal{A}(\mathbf{T},E)$.

We remark that the image of $\mathcal{B}(\mathbf{T},E)$ in
$\mathcal{B}(\mathbf{G},E)$ is the apartment
$\mathcal{A}(\mathbf{T},E)$, while the image of
$\mathcal{B}(\mathbf{T},k)$ in $\mathcal{B}(\mathbf{G},k)$ is the
set $\mathcal{A}(\mathbf{T},k)$.

\subsection{Moy-Prasad filtrations.} Regarding $\mathbf{G}$ as a
group over $E$, Moy and Prasad (see~\cite{Moy-Prasad-1}
and~\cite{Moy-Prasad-2}) define lattices in {\boldmath
$\mathfrak{g}$}$(E)$ and subgroups of $\mathbf{G}(E)$.

We can and will normalize (with respect to the normalized
valuation $\nu$) the indexing $(x,r)\in \mathcal{B}(\mathbf{G},E)
\times \mathbb{R}$ of these lattices and subgroups as in ~\cite[\S
1.2]{Fiona-Kim}. We will denote the (normalized) lattices by
{\boldmath $\mathfrak{g}$}$(E)_{x,r}$, and the (normalized)
subgroups by $\mathbf{G}(E)_{x,|r|}$.

If $\varpi_E$ is a uniformizing element of $E$, and $e=e(E/k)$ is
the ramification index of $E$ over $k$, then these normalized
lattices (resp.\ subgroups) satisfy $\varpi_E$ {\boldmath
$\mathfrak{g}$}$(E)_{x,r}=${\boldmath
$\mathfrak{g}$}$(E)_{x,r+\frac{1}{e}}$. Write $\mathfrak{g}_{x,r}$
(resp.\ $G_{x,|r|}$ ) for {\boldmath $\mathfrak{g}$}$(k)_{x,r}$
(resp.\ $\mathbf{G}(k)_{x,|r|}$).

The above normalization was chosen to have the following
property~\cite[1.4.1]{Adler}: when $E$ is a tamely ramified Galois
extension of $k$ and $x \in \mathcal{B}(\mathbf{G},k)\subset
\mathcal{B}(\mathbf{G},E)$, we have
\begin{eqnarray}\label{filtrations}
\begin{array}{cc}
\mbox{\boldmath $\mathfrak{g}$}_{x,r} = \mbox{\boldmath
$\mathfrak{g}$}(E)_{x,r} \cap \mathfrak{g}, & \mbox{ and (for } r
> 0)\;\; {G}_{x,r} = {\mathbf{G}}(E)_{x,r} \cap G.
\end{array}
\end{eqnarray}
We will also use the following notation. Let $r\in \mathbf{R}$ and
$x\in \mathcal{B}(G)$.
\begin{itemize}
\item{}$\mathfrak{g}_{x,r+} = \cup_{s > r}
\mathfrak{g}_{x,s}$ and $G_{x,|r|+} = \cup_{s > |r|} G_{x,s}$.
%\item{}$\mathfrak{g}_r=\cup_{x\in\mathcal{B}(G)}\mathfrak{g}_{x,r}$ and $\mathfrak{g}_{r+}=\cup_{s > r}\mathfrak{g}_s$.
\item{}$G_r=\cup_{x\in\mathcal{B}(G)}G_{x,r}$ and $G_{r+}=\cup_{s > r}G_s$ for $r \geq 0$.
\end{itemize}
The lattices $\mathfrak{g}_{x,r+}$ (resp.\ groups $G_{x,|r|+}$)
have analogous properties to those of $\mathfrak{g}_{x,r}$ (resp.\
$G_{x,|r|}$).
%It is proven in~\cite{Adler-DeBacker} that
%$\mathfrak{g}_r$ (resp.\ $G_{|r|}$) is a $G$-domain: a
%$G$-invariant, open and closed subset of $\mathfrak{g}$ (resp.\ $G$).
The set $G_0$ is the set of compact elements $G_{cpt}$. We
remark that $G_{cpt} \subset {\mathbf{G}}(E)_{cpt} \cap G$, and in
general they need not be equal~\cite[\S2.2.3]{Adler-DeBacker-2}.

\begin{lem}\label{fixed}
Let $\gamma$ be a compact regular semisimple element, and consider
the maximal $k$-torus $\mathbf{T}:=C_{\mathbf{G}}(\gamma)^\circ$.
Suppose that $\mathbf{T}$ splits over a tamely ramified finite
Galois extension $E$ of $k$. Then $\gamma$ fixes
$\mathcal{B}(\mathbf{T},k)$ pointwise.
\end{lem}

\begin{proof}
Recall that $\gamma$ acts on $\mathcal{A}(\mathbf{T},E)$ by
translations ~\cite[\S1]{Tits}. Since $\gamma$ belongs to a
compact subgroup, it has a fixed point $x\in
\mathcal{B}(\mathbf{G},E)$.

If $\gamma$ acts on $\mathcal{A}(\mathbf{T},E)$ by a nontrivial
translation, then for any $y\in \mathcal{A}(\mathbf{T},E)$ there
is an $n\in\mathbb{N}$ such that $d(x,y)\neq d(x,\gamma^n.y)$.
This contradicts the fact that the action preserves distances. So
$\gamma$ must act trivially on $\mathcal{A}(\mathbf{T},E)$. In
particular, $\gamma$ fixes $\mathcal{A}(\mathbf{T},k)$, and hence
$\mathcal{B}(\mathbf{T},k)$, pointwise.
\end{proof}

\subsection{Root decomposition}\label{root} Let $\mathbf{T}$ be a maximal
$k$-torus in $\mathbf{G}$ that splits over a tamely ramified
finite Galois extension $E$ of $k$. Let $\Phi(\mathbf{T},E)$
denote the set of roots of $\mathbf{G}$ with respect to $E$ and
$\mathbf{T}$, and let $\Psi(\mathbf{T},E)$ denote the
corresponding set of affine roots of $\mathbf{G}$ with respect to
$E$, $\mathbf{T}$ and $\nu$. When $\mathbf{T}$ is $k$-split, we
also write $\Phi(\mathbf{T})$ for $\Phi(\mathbf{T},k)$ (resp.\
$\Psi(\mathbf{T})$ for $\Psi(\mathbf{T},k)$). If $\psi \in
\Psi(\mathbf{T},E)$, let $\dot{\psi}\in \Phi(\mathbf{T},E)$ be the
gradient of $\psi$, and let {\boldmath
$\mathfrak{g}$}$(E)_{\dot{\psi}}\subset \mbox{\boldmath
$\mathfrak{g}$}(E)$ be the root space corresponding to
$\dot{\psi}$. We denote the root lattice in $\mbox{\boldmath
$\mathfrak{g}$}(E)_{\dot{\psi}}$ corresponding to $\psi$ by
$\mbox{\boldmath
$\mathfrak{g}$}(E)_{\psi}$~\cite[3.2]{Moy-Prasad-1}.

For $x\in\mathcal{A}(\mathbf{T},E)$ and $r\in \mathbb{R}$, let
$\mbox{\boldmath $\mathfrak{t}$}(E)_r:=\mbox{\boldmath
$\mathfrak{t}$}(E)\cap \mbox{\boldmath $\mathfrak{g}$}(E)_{x,r}$
and $\mbox{\boldmath $\mathfrak{t}$}(E)_{r+}:=\mbox{\boldmath
$\mathfrak{t}$}(E)\cap \mbox{\boldmath $\mathfrak{g}$}(E)_{x,r+}$.
Note that $\mbox{\boldmath $\mathfrak{t}$}(E)_r$ and
$\mbox{\boldmath $\mathfrak{t}$}(E)_{r+}$ are defined independent
of the choice of $x \in \mathcal{A}(\mathbf{T},E)$. Similarly one
defines the subgroups $\mathbf{T}(E)_r$ and $\mathbf{T}(E)_{r+}$
for $r \geq 0$; they have analogous properties. Note that using
our conventions we will sometimes denote $\mathbf{T}(E)_0$ by
$\mathbf{T}(E)_{cpt}$.

An alternative description is~\cite[\S2.1]{Fiona-Kim}: for $r\in
\mathbb{R}$,
$$ \mbox{\boldmath
$\mathfrak{t}$}(E)_r=\{\Gamma \in \mbox{\boldmath
$\mathfrak{t}$}(E)|\; \nu(d\chi (\Gamma))\geq r \; \mbox{ for all
}\; \chi \in \mathbf{X}^*(\mathbf{T},E) \} $$ and for $r > 0$,
$$\mathbf{T}(E)_r=\{t \in \mathbf{T}(E)|\; \nu(\chi (t)-1)\geq r
\; \mbox{ for all }\; \chi \in \mathbf{X}^*(\mathbf{T},E) \}. $$

Since $\mathbf{G}$ splits over $E$, we have
\begin{eqnarray*}
\mbox{\boldmath $\mathfrak{g}$}(E)_{x,r}&=&\mbox{\boldmath
$\mathfrak{t}$}(E)_r \oplus
\underset{\psi\in\Psi(\mathbf{T},E),\psi(x)\geq
r}\sum \mbox{\boldmath $\mathfrak{g}$}(E)_{\psi}\;,\\
\mbox{\boldmath $\mathfrak{g}$}(E)_{x,r+}&=&\mbox{\boldmath
$\mathfrak{t}$}(E)_{r+} \oplus
\underset{\psi\in\Psi(\mathbf{T},E),\psi(x)> r}\sum
\mbox{\boldmath $\mathfrak{g}$}(E)_{\psi}\;.
\end{eqnarray*}

Let $\mathfrak{t}:=\mbox{Lie}(T)$, and define
$\mathfrak{t}^\bot:=(\mbox{Ad}(\gamma)-1)\mathfrak{g}$. We have
the following decomposition~\cite[\S 18]{Queens}
\begin{eqnarray}\label{perp-decomp}
\mathfrak{g}=\mathfrak{t}\oplus\mathfrak{t}^\bot\;.
\end{eqnarray}
We write $X=Y+Z$ with respect to this decomposition; when
convenient, we also write $X_\mathfrak{t}$ for $Y$.

Fix $x\in\mathcal{B}(\mathbf{T},k)\subset
\mathcal{B}(\mathbf{G},k)$ and $r\in \mathbb{R}$. Write
$\mathfrak{t}_{r}$ for $\mathfrak{t}\cap\mathfrak{g}_{x,r}$
(resp.\ $\mathfrak{t}_{r+}$ for
$\mathfrak{t}\cap\mathfrak{g}_{x,r+}$); as mentioned earlier,
these definitions are independent of $x$. Define
$\mathfrak{t}^{\bot}_{x,r}:=\mathfrak{t}^{\bot}\cap
\mathfrak{g}_{x,r}$ (resp.\
$\mathfrak{t}^{\bot}_{x,r+}:=\mathfrak{t}^{\bot}\cap
\mathfrak{g}_{x,r+}$). We have~\cite[1.9.3]{Adler},
\begin{eqnarray}\label{perp-decomp-for-filtation}
\begin{array}{r@{\quad = \quad}l}
\mathfrak{g}_{x,r} &\mathfrak{t}_{r}\oplus
\mathfrak{t}^{\bot}_{x,r}\;,\\
\mathfrak{g}_{x,r+} & \mathfrak{t}_{r+}\oplus
\mathfrak{t}^{\bot}_{x,r+}\;.
\end{array}
\end{eqnarray}

\subsection{Hypotheses}

\noindent (HB) There is a nondegenerate $G$-invariant symmetric
bilinear form $B$ on $\mathfrak{g}$ such that we can identify
$\mathfrak{g}^*_{x,r}$ with $\mathfrak{g}_{x,r}$ via the map
$\Omega :\mathfrak{g}\rightarrow\mathfrak{g}^*$ defined by
$\Omega(X)(Y)=B(X,Y)$.\\

Groups satisfying the above hypothesis are discussed in~\cite[\S 4]{Adler-Roche}.\\

Fix $r\in\mathbb{R}_{>0}$ and $x\in\mathcal{B}(\mathbf{G},k)$. For
any $r \leq t \leq 2r$ the group $(G_{x,r}/G_{x,t})$ is abelian.
%let $(G_{x,r}/G_{x,t})^{\wedge}$ denote its Pontrjagin dual.
By hypothesis (HB), there exists a ($G_{x,0}$-equivariant)
isomorphism (see~\cite[\S 1.7]{Adler} or~\cite[p.16]{Murnaghan})
\begin{eqnarray}\label{Kirillov}
(G_{x,r}/G_{x,t})^{\wedge{\empty}} \; \cong \;
\mathfrak{g}_{x,(-t)+}/\mathfrak{g}_{x,(-r)+}\;.
\end{eqnarray}\\

\section{Regular depth}\label{regular_depth}
From now on let $\gamma \in G^{reg}$,
%(resp.\ $\lambda \in \mathfrak{g}^{reg}$)
and assume that the $k$-torus
$\mathbf{T}:=C_\mathbf{G}(\gamma)^\circ$
%(resp.\ $C_\mathbf{G}(\lambda)^\circ$)
splits over a tamely ramified finite Galois extension $E$ of $k$.
We attach to $\gamma$
%(resp.\ $\lambda$)
the following rational number $s(\gamma)$.
\begin{defn}\label{s}
For each $\alpha\in \Phi(\mathbf{T},E)$ let
$s_\alpha(\gamma):=\nu(\alpha(\gamma)-1)$ and
%and let $s_\alpha(\lambda):=\nu(d\alpha(\lambda))$ ($d\alpha$ is the derivative of $\alpha$).
define $s(\gamma):=\mbox{max}\{\,s_\alpha(\gamma)\,|\,\alpha\in
\Phi(\mathbf{T},E)\,\}$.
%and $s(\lambda):=\mbox{max}\{\,s_\alpha(\lambda)\,|\,\alpha\in \Phi(\mathbf{T},E)\,\}$.
\end{defn}

\begin{rmk}
Note that $s(\gamma )$ is not the same as the depth of $\gamma$
(as defined in~\cite{Adler-DeBacker}). But for good
elements~\cite[\S2.2]{Adler}, these two notions agree.
\end{rmk}

\begin{rmk}\label{regrmk}
A priori $s( \gamma ) \in \mathbb{Q}\cup \{+ \infty \}$, but since
$\gamma$ is regular, $\alpha(\gamma)\neq 1$ for all $\alpha\in
\Phi(\mathbf{T},E)$ and so $s(\gamma)\in \mathbb{Q}$. If $\gamma$
is compact then $s(\gamma)\geq 0$. Also note that $s(\gamma
z)=s(\gamma)$ for all $z$ in the center $Z(G)$ of $G$ and that
$s(g \gamma g^{-1}) = s(\gamma )$ for all $g\in G$.
\end{rmk}

We will need the following basic properties of $s(\gamma)$.

\begin{lem}\label{deepness}
Suppose $\gamma\in T_{cpt}$ and $\gamma'\in T_{s(\gamma)+}$.
\begin{enumerate}
\item{} $s(\gamma\gamma')=s(\gamma)$ and for  $\alpha\in\Phi(\mathbf{T},E)$,
we have $|\alpha(\gamma \gamma ' )-1|=|\alpha(\gamma )-1|$.
\item{} $\gamma\gamma' \in T_{cpt}$.
\end{enumerate}
\end{lem}

\begin{proof}
\begin{enumerate}
\item{}
Fix $r > s(\gamma) \geq 0$ such that $T_r = T_{s(\gamma)+}$. With
this notation $\gamma'\in T_r$. By the alternative description of
$T_r$, for any $\chi \in \mathbf{X}^*(\mathbf{T},E)$, $\chi(\gamma
')=1+\mu '$ where $\nu (\mu ')\geq r$. Thus for any $\alpha \in
\Phi(\mathbf{T},E)$, $\alpha(\gamma ')=1+\lambda '$ where $\nu
(\lambda ')\geq r$.

Note that since each $\alpha \in \Phi(\mathbf{T},E)$ is
continuous, $\alpha (T(E)_{cpt})\subset R_E^\times$. Since $\gamma
\in T_{cpt} \subset T(E)_{cpt}$ we get that $\alpha (\gamma)$ is a
unit.

Now $\alpha (\gamma\gamma ')-1=\alpha (\gamma)\alpha (\gamma
')-1=\alpha (\gamma )(1+\lambda ')-1= (\alpha (\gamma)-1)+\alpha
(\gamma)\lambda '$. Using $\nu (\alpha (\gamma)-1)=:s_\alpha
(\gamma)$, $\alpha (\gamma)$ is a unit, and $\nu (\lambda ') \geq
r > s(\gamma) \geq s_\alpha (\gamma)$, we have $\nu (\alpha
(\gamma\gamma ')-1)=\nu (\alpha (\gamma)-1))$ (or equivalently
$|\alpha (\gamma\gamma ')-1|=|\alpha (\gamma)-1)|$) for all
$\alpha\in\Phi(\mathbf{T},E)$. Thus $s(\gamma\gamma
'):=\underset{\alpha}{\mbox{max}}\{\nu (\alpha (\gamma\gamma
')-1)\}= \underset{\alpha}{\mbox{max}}\{\nu (\alpha
(\gamma)-1)\}=: s (\gamma)$.

\item{} Since $\gamma$ and $\gamma '$ are in $T_{cpt}$, so is
their product.
%Jeff: Since $s(\gamma ) \geq \mbox{depth}(\gamma)$, $\gamma '
%\in T_{s(\gamma )+} \subseteq  T_{\mbox{depth}(\gamma )+}$. Hence
%$\gamma$ and $\gamma\gamma '$ have the same depth; in particular
%$\gamma \in G_{cpt}$ if and only if $\gamma\gamma ' \in G_{cpt}$.
\end{enumerate}
\end{proof}

\begin{cor}\label{deepness cor}
Let $\gamma \in T$ be a compact regular semisimple element. Then
$\gamma T_{s(\gamma)+}\subset G^{reg}$.
\end{cor}
\begin{proof}
For $t \in T \cap G^{reg}$, following~\cite[\S18]{Queens}, define
$$D_{G/T}(t):=\mbox{det}\left(\mbox{Ad}(t)-1)\right)|_{\mathfrak{g}/\mathfrak{t}}=\underset{\alpha\in\Phi
(\mathbf{T},E) }\prod (\alpha (t)-1).$$ Then $t\in T\cap
G^{reg}\Leftrightarrow D_{G/T}(t)\neq 0 \Leftrightarrow
|D_{G/T}(t)|\neq 0$. Using Lemma~\ref{deepness} with $\gamma \in
T\cap G_{cpt}$ and $\gamma '\in T_{s(\gamma)+}$, we get
$|D_{G/T}(\gamma\gamma ' )|=\underset{\alpha}\prod |\alpha
(\gamma\gamma ')-1|=\underset{\alpha}\prod |\alpha
(\gamma)-1|=|D_{G/T}(\gamma )|\neq 0$.
\end{proof}

\section{Some Technical Lemmas}\label{tech_lemma}
The next lemma will generalize the following example.

\begin{ex}
$\mathbf{G}=\mathbf{GL_2}$, $\mathbf{T}$ a $k$-split maximal
torus. Choose $x_0\in\mathcal{B}(\mathbf{G},k)$ so that
$G_{x_0,0}=GL_2(R)$. Any $X\in
\mathcal{N}\cap(\mathfrak{g}_{x_0,r}\smallsetminus
\mathfrak{g}_{x_0,r+})$ is of the form ${\empty}^k
\begin{pmatrix} 0 & x \\ 0 & 0 \end{pmatrix}$, for some $k\in G_{x_0,0}$
(see~\cite[9.2.1]{Singapore}). Thus
\begin{eqnarray*}
X &=& \begin{pmatrix} a & b \\ c & d \end{pmatrix}
\begin{pmatrix} 0 & x \\ 0 & 0 \end{pmatrix}
\begin{pmatrix} a & b \\ c & d \end{pmatrix}^{-1}\\
&=& \frac{x}{ad-bc}
\begin{pmatrix} -ac & a^2 \\ -c^2 & ac \end{pmatrix}.
\end{eqnarray*}
Write $X=Y+Z$ as in ~(\ref{perp-decomp}) and note that the depth
of $X$ with respect to the filtration
$\{\mathfrak{g}_{x_0,r}\}_{r\in\mathbb{R}}$ of $\mathfrak{g}$ is
controlled by $Z$. This is the case since max$\{
\nu(a^2),\nu(-c^2)\}\,\geq\, \nu(ac)$ and $ad-bc\in R^\times$.
\end{ex}

\begin{lem}\label{keylemma}
Fix $x\in \mathcal{B}(\mathbf{T},k)$ and $r\in \mathbb{R}$. For
$X\in\mathcal{N}\cap(\mathfrak{g}_{x,r}\smallsetminus
\mathfrak{g}_{x,r+})$, write $X=Y+Z$ as in ~(\ref{perp-decomp}).
Then $Z\in \mathfrak{g}_{x,r}\smallsetminus \mathfrak{g}_{x,r+}$.
\end{lem}

\begin{proof}
We first prove the case when the maximal $k$-torus $\mathbf{T}$ is
$k$-split and then reduce the general case to this case.

{\bf Split case.} Assume $\mathbf{T}$ is $k$-split.
%so $\mathbf{G}$ regarded as group over $k$ is a split group
Note that $\mathfrak{t}^\bot=\oplus_{\alpha\in\Phi(\mathbf{T})}
\mathfrak{g}_\alpha$. Fix a system of simple roots $\Delta$ in
$\Phi(\mathbf{T})$ and choose a Chevalley basis for $\mathfrak{g}$
as in~\cite[\S 1.2]{Adler}. Such a basis contains elements $H_b$
and $E_b$ in $\mathfrak{g}$ for each $b\in \Phi(\mathbf{T})$. If
$\mathbf{G}$ is semisimple, then the set $\{H_b|\;b \in \Delta
\}\cup \{E_b|\;b\in\Phi(\mathbf{T}) \}$ is a basis for
$\mathfrak{g}$. These elements also satisfy the commutation
relations listed in~\cite[1.2.1]{Adler}. With respect to this
choice of Chevalley basis, the adjoint representation is
determined by the following formulas~\cite[1.2.5]{Adler}:

\begin{eqnarray}\label{Chevalley}
\left\{
\begin{array}{r@{\quad = \quad}l}
\mbox{Ad}(e_b(\lambda))E_c & \left\{ \begin{array}{l}
                             E_b\\
                             E_c+\lambda H_b - \lambda^2E_b\\
                             \sum_{i\geq0}M_{b,c;i}\,\lambda^i\,
                            E_{ib+c}\\
                            \end{array}
                            \begin{array}
                            {r@{\quad if \quad}l}
                            & c=b\\
                            & c=-b\\
                            & c \neq \pm b
                            \end{array}
                       \right.\\
\mbox{Ad}(t)E_c & c(t)E_c\\
\mbox{Ad}(e_b(\lambda))H & H - db(H)\lambda E_b\\
\mbox{Ad}(t)H & H
\end{array}
\right.
\end{eqnarray}
for all $H\in$ Lie$(T)$, all $t\in T$ and all $\lambda\in k$. Here
$e_b$ is the unique map $e_b:
\mathbf{Add}\longrightarrow\mathbf{G}$ such that d$e_b(1)=E_b$
($de_b$ is the derivative of $e_b$); and $M_{b,c;i}$ are constants
with $M_{b,c;0}=1$.

Let $B$ be the Borel subgroup associated to $\Delta$ (with Levi
decomposition $B=TN$ and opposite Borel
$\overline{B}=T\overline{N}$). We have
$\mathfrak{g}=\mathfrak{t}\oplus\mathfrak{n}\oplus\overline{\mathfrak{n}}$,
where $\mathfrak{n}:=\mbox{Lie}(N)$ and
$\overline{\mathfrak{n}}:=\mbox{Lie}(\overline{N})$. Note that
$\mathfrak{n}\oplus\overline{\mathfrak{n}}=\oplus_{\alpha\in\Phi(\mathbf{T})}
\mathfrak{g}_\alpha = \mathfrak{t}^\bot$. Recall that $G_{x,0}$
acts on $\mathfrak{g}_{x,r}$ (and on $\mathfrak{g}_{x,r+}$).

Given $X\in\mathcal{N}\cap(\mathfrak{g}_{x,r}\smallsetminus
\mathfrak{g}_{x,r+})$, we can use
~\cite[Proposition~3.5.1]{Adler-DeBacker} (with $T$ playing the
role of $M$) to conclude that there exists a group element $n\in N
\cap G_{x,0}$ such that $(^nX)_\mathfrak{t}\in \mathfrak{t}_{r+}$
(where $^nX$ denotes $\mbox{Ad}(n)X$).

Write $X=Y+Z$ as in ~(\ref{perp-decomp}) and assume for a
contradiction that $Z\in\mathfrak{g}_{x,r+}$. Since $X\in
\mathfrak{g}_{x,r}\smallsetminus \mathfrak{g}_{x,r+}$, the
assumption implies that $Y\in \mathfrak{t}\cap
(\mathfrak{g}_{x,r}\smallsetminus
\mathfrak{g}_{x,r+})=\mathfrak{t}_{r}\smallsetminus
\mathfrak{t}_{r+}$.

Using the properties (\ref{Chevalley}) of the Chevalley basis, one
can easily check that the set
$(\mathfrak{t}_r\smallsetminus\mathfrak{t}_{r+})\oplus\mathfrak{n}$
is preserved under the action of $\mbox{Ad}(e_b(\lambda))$ for all
$b\in\Phi^+(\mathbf{T})$, where $\Phi^+(\mathbf{T})$ are the
positive roots with respect to $\Delta$. Since $\{e_b(\lambda)
\;|\;b\in \Phi^+(\mathbf{T}) \}$ generates $N$, we conclude that
$^nY\in
(\mathfrak{t}_r\smallsetminus\mathfrak{t}_{r+})\oplus\mathfrak{n}$,
 and hence that $(^nY)_\mathfrak{t} \in
\mathfrak{t}_r\smallsetminus\mathfrak{t}_{r+}$.

On the other hand we have ${^nX}={^nY}+{^nZ}$, where
$^nZ\in\mathfrak{g}_{x,r+}$. Taking the $\mathfrak{t}$ components,
we get, $(^nX)_\mathfrak{t}=(^nY)_\mathfrak{t} +
(^nZ)_\mathfrak{t}$, with $(^nZ)_\mathfrak{t}\in
\mathfrak{t}_{r+}$. Since $(^nX)_\mathfrak{t}\in
\mathfrak{t}_{r+}$, we conclude that $(^nY)_\mathfrak{t}\in
\mathfrak{t}_{r+}$. This contradicts $(^nY)_\mathfrak{t}\in
\mathfrak{t}_{r}\smallsetminus \mathfrak{t}_{r+}$.

Hence $Z\in \mathfrak{g}_{x,r}\smallsetminus \mathfrak{g}_{x,r+}$
(note that from the decomposition
(\ref{perp-decomp-for-filtation}) it is clear that $Z\in
\mathfrak{g}_{x,r}$).

{\bf General case.} We now assume $\mathbf{T}$ is an $E$-split
maximal $k$-torus. Define $\qquad \mathfrak{t}(E)^\bot
:=(\mbox{Ad}(\gamma)-1)\mathfrak{g}(E)$. We have the following
analogue of (\ref{perp-decomp})
\begin{eqnarray}\label{E-perp-decomp}
\mbox{\boldmath $\mathfrak{g}$}(E)=\mbox{\boldmath
$\mathfrak{t}$}(E)\oplus\mbox{\boldmath $\mathfrak{t}$}(E)^\bot.
\end{eqnarray}
Note that $\mathfrak{t}\subset \mbox{\boldmath $\mathfrak{t}$}(E)$
and $\mathfrak{t}^\bot \subset \mbox{\boldmath
$\mathfrak{t}$}(E)^\bot$. So the decomposition $X=Y+Z$ (as in
~(\ref{perp-decomp})) for $X\in\mathfrak{g}$ is the same whether
viewed in $\mathfrak{g}$ or in $\mbox{\boldmath
$\mathfrak{g}$}(E)$.

Since $X\in \mathfrak{g}_{x,r}\smallsetminus \mathfrak{g}_{x,r+}$,
equations (\ref{filtrations}) imply that $X\in (\mbox{\boldmath $
\mathfrak{g}$}(E)_{x,r}\smallsetminus \mbox{\boldmath $
\mathfrak{g}$}(E)_{x,r+})\cap\mathfrak{g}$. Since
$X\in\mathcal{N}\subset\mbox{\boldmath $\mathcal{N}$}(E)$ (where
$\mbox{\boldmath $\mathcal{N}$}(E)$ is the set of nilpotent
elements in $\mbox{\boldmath $\mathfrak{g}$}(E))$, we have that
$X\in \mbox{\boldmath $\mathcal{N}$}(E)\cap(\mbox{\boldmath
$\mathfrak{g}$}(E)_{x,r}\smallsetminus \mbox{\boldmath
$\mathfrak{g}$}(E)_{x,r+})$. Now since $\mathbf{T}$ splits over
$E$ we can regard $\mathbf{G}$ over $E$ as a split group and hence
apply all the constructions of the split case above. So by the
considerations of the split case above we conclude that
$Z\in\mbox{\boldmath$\mathfrak{g}$}(E)_{x,r}\smallsetminus
\mbox{\boldmath $\mathfrak{g}$}(E)_{x,r+}$. Intersecting with
$\mathfrak{g}$ gives $Z\in\mathfrak{g}_{x,r}\smallsetminus
\mathfrak{g}_{x,r+}$.
\end{proof}

From now on we assume that $\gamma$ is also compact.
%$\gamma\in G_0$.
Recall that this implies that $s(\gamma)\geq 0$ (see
Remark~\ref{regrmk}).

\begin{lem}\label{Z-gamma} Let $t\in \mathbb{R}$ and $x\in\mathcal{B}(\mathbf{T},k)$.
If $Z\in \mathfrak{t}^\bot\cap (\mathfrak{g}_{x,-t}\smallsetminus
\mathfrak{g}_{x,(-t)+})$ then ${^\gamma Z}-Z \not\in
\mathfrak{g}_{x,(-t+s(\gamma))+}$.
\end{lem}

\begin{proof}
Using the root decomposition $\mbox{\boldmath
$\mathfrak{t}$}(E)^\bot=\oplus_{\alpha\in\Phi(\mathbf{T},E)}
\mbox{\boldmath $\mathfrak{g}$}(E)_\alpha$, for $Z\in
\mathfrak{t}^\bot \subset \mbox{\boldmath $\mathfrak{t}$}(E)^\bot$
we write $Z=\sum Z_\alpha$. Then ${}^\gamma Z-Z = \sum ({}^\gamma
Z_\alpha - Z_\alpha ) = \sum(\alpha (\gamma)-1)Z_\alpha $.

By assumption $Z\not\in \mathfrak{g}_{x,(-t)+}$, hence (see
equations (\ref{filtrations})) $Z\not\in \mbox{\boldmath
$\mathfrak{g}$}(E)_{x,(-t)+}$. Thus for some $\alpha \in
\Phi(\mathbf{T},E)$, $Z_\alpha\not\in \mbox{\boldmath
$\mathfrak{g}$}(E)_{x,(-t)+}$, and so by definition of
$s_\alpha(\gamma)$, $(\alpha(\gamma)-1)Z_\alpha\not\in
\mbox{\boldmath $\mathfrak{g}$}(E)_{x,(-t+s_\alpha(\gamma))+}$. It
follows by definition of $s(\gamma)$, that
$(\alpha(\gamma)-1)Z_\alpha\not\in \mbox{\boldmath $
\mathfrak{g}$}(E)_{x,(-t+s(\gamma))+}$. Hence ${}^\gamma Z-Z =
\sum (\alpha (\gamma)-1) Z_\alpha \not\in \mbox{\boldmath $
\mathfrak{g}$}(E)_{x,(-t+s(\gamma))+}$. Intersecting with
$\mathfrak{g}$ we get that ${}^\gamma Z-Z \not\in
\mathfrak{g}_{x,(-t+s(\gamma))+}$.
\end{proof}

\begin{propn}\label{proposition}
Let $r\in \mathbb{R}$ and $x\in\mathcal{B}(\mathbf{T},k)$. If
$X\in\mathcal{N}\cap \mathfrak{g}_{x,(-2r)+}$ satisfies ${^\gamma
X}-X \in \mathfrak{g}_{x,(-r)+}$, then $X\in
\mathfrak{g}_{x,(-r-s(\gamma))+}$.
\end{propn}

\begin{proof}
Fix $t < 2r$ such that $X\in\mathcal{N}\cap
(\mathfrak{g}_{x,-t}\smallsetminus \mathfrak{g}_{x,(-t)+})$.

Write $X=Y+Z$ as in ~(\ref{perp-decomp}). By Lemma~\ref{keylemma},
$Z\in\mathfrak{t}^\bot\cap (\mathfrak{g}_{x,-t}\smallsetminus
\mathfrak{g}_{x,(-t)+})$, and so by Lemma~\ref{Z-gamma}, ${^\gamma
Z}-Z \not\in \mathfrak{g}_{x,(-t+s(\gamma))+}$.

On the other hand, since $\gamma$ acts trivially on $Y$ (because
$Y\in \mathfrak{t}=C_{\mathfrak{g}}(\gamma)$), ${^\gamma Z}-Z =
{^\gamma X}-X \in \mathfrak{g}_{x,(-r)+}$.

Thus $-t+s(\gamma) > -r$, or equivalently $-t > -r-s(\gamma)$,
which implies that $X\in\mathfrak{g}_{x,-t}\subseteq
\mathfrak{g}_{x,(-r-s(\gamma))+}$.
\end{proof}

\begin{defn}
A character $d\in (G_{x,r}/G_{x,2r})^{\wedge}$ is called
degenerate if under the isomorphism (\ref{Kirillov}), the
corresponding coset $X+\mathfrak{g}_{x,(-r)+}$ contains nilpotent
elements.
\end{defn}

\begin{defn} Let $K$ be a compact subgroup of $G$ and $d\in K^{\wedge}$. For $g\in G$, let ${}^g
d$ denote the representation of $g K g ^{-1}$ defined as ${}^g d
(g k g ^{-1}):=d(k)$. We say that $g$ intertwines $d$ with itself
if upon restriction to $g K g ^{-1} \cap K$, $d$ and ${}^g d$
contain a common representation (up to isomorphism) of $g K g
^{-1} \cap K$.
\end{defn}

\begin{cor}\label{cor to propn}
Let $x\in\mathcal{B}(\mathbf{T},k)$, $r\in \mathbb{R}_{> 0}$, and
assume $d\in (G_{x,r}/G_{x,2r})^{\wedge}$ is degenerate. If
$\gamma$ intertwines $d$ with itself then $\;d\in
(G_{x,r}/G_{x,r+s(\gamma)})^{\wedge}$.
\end{cor}

\begin{proof}
Let $X+\mathfrak{g}_{x,(-r)+}$ be the coset in
$\mathfrak{g}_{x,(-2r)+}/\mathfrak{g}_{x,(-r)+}$ corresponding to
$d$ under the isomorphism (\ref{Kirillov}). Since this coset is
degenerate, we can assume that $X\in \mathcal{N}$.

Since $\gamma$ fixes $x$ (Lemma~\ref{fixed}), $\gamma$ stabilizes
$G_{x,r}$.
% and so $G_{x,r}\cap \gamma G_{x,r} \gamma^{-1}=G_{x,r}$
Thus having $\gamma$ intertwine $d$ with itself amounts to having
$d \cong {}^\gamma d$; or furthermore, since $d$ is
one-dimensional, $d = {}^\gamma d$. Under the isomorphism
(\ref{Kirillov}), we get $X+\mathfrak{g}_{x,(-r)+} = {}^\gamma
(X+\mathfrak{g}_{x,(-r)+})$, or equivalently that ${}^\gamma X -
X\in \mathfrak{g}_{x,(-r)+}$. Now apply
Proposition~\ref{proposition} to conclude that $X\in
\mathfrak{g}_{x,(-r-s(\gamma))+}$, which under the isomorphism
(\ref{Kirillov}) gives that $d\in
(G_{x,r}/G_{x,r+s(\gamma)})^{\wedge}$.
\end{proof}

\section{Partial Traces}\label{partial traces}
Let $(\pi, V)$ be an irreducible admissible representation of $G$.
Let $K$ be an open compact subgroup of $G$. Let $V=\bigoplus_{d\in
K^\wedge} V_d$ be the decomposition of $V$ into $K$-isotypic
components. Let $E_d$ denote the $K$-equivariant projection from
$V$ to $V_d$. For $f\in C_c^\infty (G)$ define the distribution
$\Theta_d (f):= \mathrm{trace}\,(E_d \,\pi (f)\, E_d)$, the
`partial trace of $\pi$ with respect to $d$'. The distribution
$\Theta_d$ is represented by the locally constant function
$\Theta_d (x):= \mathrm{trace}\, (E_d\, \pi (x) \,E_d)$ on $G$.
Recall that it is known that the distribution $\Theta_\pi (f):=
\mathrm{trace}\, \pi (f)$ is also represented by a locally
constant function, $\Theta_\pi$, on $G^{reg}$; we will not use
this fact here. It follows from the definitions that as
distributions
$$\Theta_\pi (f)=\sum_{d\in K^\wedge}\Theta_d (f) \;\mbox{ for all } f\in
C_c^\infty (G).$$

\begin{rmk}\label{supp_of d}
For (some) $\omega \subset G^{reg}$ compact, Corollary~\ref{cordd}
and the proof of Theorem~\ref{mainresult} will imply that, for all
$f \in C_c^{\infty}(\omega)$, this sum is {\it finite}.
\end{rmk}

\begin{lem}\label{conjlem} $\Theta_d(kxk^{-1})=\Theta_d(x)$ for all $x \in G$
and all $k \in K$.
\end{lem}
\begin{proof}
Since $E_d$ is $K$-equivariant, it commutes with $\pi(k)$ for all
$k\in K$.
\begin{eqnarray*}
\Theta_d(kxk^{-1}) &=& \mathrm{trace}(E_d \pi(kxk^{-1}) E_d)\\
                   &=& \mathrm{trace}(E_d \pi(k)\pi(x)\pi(k^{-1}) E_d)\\
                   &=& \mathrm{trace}(\pi(k)E_d\pi(x)E_d\pi(k^{-1}))\\
                   &=& \mathrm{trace}(E_d\pi(x)E_d) = \Theta_d(x).
\end{eqnarray*}
\end{proof}

Let $N$ be an open compact subgroup of $G$ which is normal in $K$.
Suppose $g \in G$ normalizes $K$ and $N$. Let $d\in K^\wedge$.
Considered as a representation of $N$, $d$ decomposes into a
finite sum of irreducible representations
$$d_1\oplus\cdots\oplus d_n\,.$$

\begin{propn}\label{intertwine} Suppose $\Theta_d (g) \neq 0$. Then $d \cong {}^g
d$ as representations of $K$ and also for some
$i\in\{1,\cdots,n\}$, $d_i \cong {}^g d_i$ as representations of
$N$.
\end{propn}

\begin{proof} We refer to the appendix.
Since $g$ permutes the $V_{d'}$'s (Theorem~\ref{Clifford}.1),
$0\neq \Theta_d (g)=\mathrm{trace} (E_d \pi (g) E_d)$ implies that
$g$ must stabilize $V_d$. Fix a decomposition $(\ddag)$ as in
Theorem~\ref{Clifford}.2, and let $E_{W_i}$ denote the
$K$-equivariant projections onto $W_i$. Since $E_d = \sum
E_{W_i}$, $\mathrm{trace} (E_d \pi (g) E_d)\neq 0$ implies that
for some $i$, $g$ must stabilize $W_i$, and that $\mathrm{trace}
(E_{W_i} \pi (g) E_{W_i})\neq 0$. By Theorem~\ref{Clifford}.3, $ d
\cong {}^g d$, which proves the first part of the theorem.

Now as a representation of $N$,
$$W_i \underset{N} = \bigoplus_{j} W_{i,d_j},$$ where $W_{i,d_j}$ are
the $d_j$-isotypic components of $W_i$. Since $g$ stabilizes $N$,
it must permute the $W_{i,d_j}$'s (Theorem~\ref{Clifford}.1).
Since $E_{W_i} = \sum E_{W_{i,d_j}}$, having $\mathrm{trace}
(E_{W_i} \pi (g) E_{W_i})\neq 0$ implies that for some $j$, $g$
must stabilize $W_{i,d_j}$, and that $\mathrm{trace}
(E_{W_{i,d_j}} \pi (g) E_{W_{i,d_j}})\neq 0$. Fix a decomposition
$(\ddag)$ as in Theorem~\ref{Clifford}.2 for $W_{i,d_j}$:
$$W_{i,d_j}\underset{N}\cong \bigoplus d_j .$$ Since $E_{W_{i,d_j}} = \sum
E_{d_j}$, $\mathrm{trace} (E_{W_{i,d_j}} \pi (g)
E_{W_{i,d_j}})\neq 0$ implies that $g$ must stabilize one of the
$d_j$'s. By Theorem~\ref{Clifford}.3, $ d_j \cong {}^g d_j$, which
proves the second part of the theorem.
\end{proof}

The following theorem and corollaries are used in the proof of
Theorem~\ref{mainresult} to show that for $f$ with compact
support, the sum $\sum_{d\in K^\wedge}\Theta_d (f)$ is finite (see
also Remark~\ref{supp_of d}).

\begin{thm}\label{Kirillov thm} Fix $x\in\mathcal{B}(\mathbf{T},k)$
and let $r > {\rm max}\{ s(\gamma),\rho (\pi) \}$. If $d \in
(G_{x,r})^{\wedge}$ satisfies $\Theta_d(\gamma)\neq 0$, then $d
\in(G_{x,r}/G_{x,r+s(\gamma)})^{\wedge}.$
\end{thm}

\begin{proof}
If $d$ is trivial we are done, so assume it is not. Let $t$ be the
smallest number such that $d|_{G_{x,t+}}$ is trivial (so in
particular $d|_{G_{x,t}}$ is nontrivial).

\textbf{Case $t < 2r$}: Pick $s \leq 2r$ such that
$G_{x,s}=G_{x,t+}$. Consider $d$ as an element of
$(G_{x,r}/G_{x,2r})^{\wedge}$. By Proposition~\ref{intertwine},
$\Theta_d(\gamma)\neq 0$ implies that $d \cong {}^\gamma d$. Also,
$\Theta_d(\gamma)\neq 0$ implies that $d \subset \pi |_{G_{x,r}}$;
since $r > \rho (\pi)$ this means that $d$ is degenerate
(see~\cite[\S7.6]{Singapore}). Now apply Corollary~\ref{cor to
propn}.

\textbf{Case $t \geq 2r $}: Note that $\frac{t}{2} \geq r >
s(\gamma)$. For $\epsilon > 0$ such that $\frac{t}{2} >
\frac{\epsilon}{2} + s(\gamma)$, let $s=t+\epsilon$. By making
$\epsilon$ smaller if necessary, we can make sure that
$G_{x,s}=G_{x,t+}$. Note that $t
> \frac{t}{2} + \frac{\epsilon}{2} + s(\gamma)=
\frac{s}{2}+s(\gamma)$.

Since $\frac{s}{2} > \frac{t}{2} \geq r$ it makes sense to
restrict $d$ to $G_{x,\frac{s}{2}}$ and think of it as an element
of $(G_{x,\frac{s}{2}}/G_{x,s})^{\wedge}$. As a representation of
$G_{x,\frac{s}{2}}/G_{x,s}$, $d$ decomposes into a finite sum of
irreducible (one-dimensional) representations
$$d_1\oplus\cdots\oplus d_n\,.$$

Let $X_i + \mathfrak{g}_{x,(-\frac{s}{2})+}$ be the coset in
$\mathfrak{g}_{x,(-s)+}/\mathfrak{g}_{x,(-\frac{s}{2})+}$
corresponding to $d_i$ under the isomorphism~(\ref{Kirillov}).

%Roger Howe calls such $d_i$ `shallow'~\cite{Howe}.

By Proposition~\ref{intertwine}, $0 \neq \Theta_d(\gamma)$ implies
that for some $j$, $d_j \cong {}^\gamma d_j$.

Now $d\subset \pi |_{G_{x,r}}$, implies that $d_j\subset \pi
|_{G_{x,\frac{s}{2}}}$ and since $\frac{s}{2} > r > \rho(\pi)$ we
have that $d_j$ is degenerate. Apply Corollary~\ref{cor to propn}
to $d_j$
%thought of as an element of$(G_{x,\frac{s}{2}}/G_{x,s})^{\wedge}$,
to conclude that
$d_j\in(G_{x,\frac{s}{2}}/G_{x,\frac{s}{2}+s(\gamma)})^{\wedge}$.
In particular $d_j$ is trivial on $G_{x,\frac{s}{2}+s(\gamma)}$,
and hence on $G_{x,t}$.

Since $G_{x,r}$ normalizes $G_{x,\frac{s}{2}}$, it acts by
permutations on the $d_i$'s. Since $d$ is irreducible, this action
is transitive. Hence all the $d_i$'s are conjugate by elements of
$G_{x,r}$. By the conjugation of the $d_i$'s and the fact that
$d_j|_{G_{x,t}}=1$ it follows that $d_i|_{G_{x,t}}=1$ for all $i$,
and so $d$ itself is trivial on $G_{x,t}$. This contradicts the
definition of $t$. Hence this case is not possible and $t < 2r$.
\end{proof}

\begin{cor}\label{cord} Fix $x\in\mathcal{B}(\mathbf{T},k)$
and let $r > {\rm max}\{ s(\gamma),\rho (\pi) \}$. Let $X$ denote
$\gamma T_{r+s(\gamma)}$, a compact subset of $T \cap G^{reg}$. If
$d \in (G_{x,r})^{\wedge}$ satisfies $\Theta_d(\gamma ')\neq 0$
for some $\gamma ' \in X$, then $d
\in(G_{x,r}/G_{x,r+s(\gamma)})^{\wedge}.$
\end{cor}

\begin{proof} Lemma~\ref{deepness} implies that $\gamma ' $ fixes
$x$ and that $s(\gamma ')$=$s(\gamma )$. Now apply
Theorem~\ref{Kirillov thm} to $\gamma '$.
\end{proof}

\begin{cor}\label{cordd} Fix $x\in\mathcal{B}(\mathbf{T},k)$
and let $r > {\rm max}\{ s(\gamma),\rho (\pi) \}$. Let $\omega$
denote ${}^{G_{x,r}}(\gamma T_{r+s(\gamma)})$, an open compact
subset of $G^{reg}$. Then $\Theta_d$ vanishes on $\omega$ for all
$d \notin(G_{x,r}/G_{x,r+s(\gamma)})^{\wedge}$. Furthermore,
$\Theta_d (x)=\Theta_d (\gamma) $ for all $x\in \omega$ and all $d
\in(G_{x,r}/G_{x,r+s(\gamma)})^{\wedge}$.
\end{cor}

\begin{proof} Follows immediately from Lemma~\ref{conjlem} and
Corollary~\ref{cord}.
\end{proof}

\section{Proof of the Main Theorem}
Let $r > {\rm max}\{ s(\gamma),\rho (\pi) \}$. Denote the finite
set $(G_{x,r}/G_{x,r+s(\gamma)})^{\wedge}$ by $F$.

\begin{thm}\label{mainresult} The distribution $\Theta_\pi$ is represented on the set
${}^{G}(\gamma T_{r+s(\gamma)})$ by a constant function.
\end{thm}

\begin{proof}
Using Corollary~\ref{cordd}, we have for all $f \in
C_c^{\infty}(G)$ whose support is contained in $\omega$,
\begin{align*}
\Theta_\pi(f) = \sum_{d\in (G_{x,r})^\wedge}\Theta_d(f)
              = \sum_{d\in F}\int_{\omega}\Theta_d(x)f(x)dx
              &= \sum_{d\in
              F}\int_{\omega}\Theta_d(\gamma)f(x)dx\\
              &= \int_{\omega}\left(\sum_{d\in
              F}\Theta_d(\gamma)\right)f(x)dx.
\end{align*}
Thus $\Theta_\pi$ is represented by the constant function
$\sum_{d\in F}\Theta_d(\gamma)$ on $\omega$, i.e.
$\Theta_\pi(x)=\sum_{d\in F}\Theta_d(\gamma)$ for all $x\in
\omega$. Since $\Theta_\pi$ is conjugation invariant, we get
$\Theta_\pi(gxg^{-1})=\Theta_\pi(x)= \sum_{d\in
F}\Theta_d(\gamma)$ for all $x\in \omega$ and all $g\in G$.
\end{proof}
\begin{rmk} This gives a new proof of the local constancy
(near compact regular semisimple tame elements $\gamma$) of the
character of an irreducible admissible representation for an
arbitrary reductive $p$-adic group $G$.
\end{rmk}

\section{Appendix}
We prove some variations of Clifford theory~\cite[\S14]{Clifford}.
Let $K$ and $N$ be open compact subgroups of $G$, such that $N$ is
a normal subgroup of $K$. Let $(\pi,V)$ be an irreducible
admissible representation of $G$ and let
\begin{gather}V\underset{K}=\bigoplus_{d\in K^\wedge}V_d \tag{$\dagger$}
\end{gather}
be the (canonical) decomposition of $V$ into $K$-isotypic
components. Here $V_d$ denotes the $d$-isotypic component of $V$,
i.e. the sum of all the $K$-submodules of $V$ isomorphic to
$d=(d,W)$. Each isotypic component $V_d$ decomposes
(non-canonically) into a finite sum of isomorphic copies $W_i
\underset{K}\cong W$ of $(d,W)$
\begin{gather} V_d\underset{K}\cong\bigoplus_{i}  W_i \;. \tag{$\ddagger$}
\end{gather}

\begin{thm}\label{Clifford}
Suppose $g \in G$ normalizes $K$ and $N$. Then
\begin{enumerate}
\item{} The action of $g$ permutes the $V_d$'s.
\item{} Suppose $g$ stabilizes $V_d$. Then there exists a
decomposition ($\ddagger$) such that the action of $g$ permutes
the $W_i$.
\item{} Suppose $W'$ is a $K$-submodule of $V$, isomorphic to $W$, and stable under the action
$g$. Then ${}^g d \cong d$.
\end{enumerate}
\end{thm}
\begin{proof}
\begin{enumerate}
\item{} This follows from the fact that for any two
$K$-submodules $W'$ and $W''$ of $V$, if $W' \underset{K}\cong
W''$ then $g W' \underset{K}\cong g W''$.
\item{} Let $W'$ be an irreducible $K$-submodule of $V_d$,
isomorphic to $W$. Since $g$ normalizes $K$ and stabilizes $V_d$,
$g W'$ is a $K$-submodule of $V_d$. Since $W'$ is irreducible, so
is $g W'$. As an irreducible submodule of $V_d$, $g W'$ must be
isomorphic to $W$. By irreducibility either $W' \cap g W' = \{0\}$
or $W' = g W'$. Thus the orbit of $W'$ under $g$ is a collection
of subspaces with trivial pairwise intersection, and so $g$ acts
on their sum as a desired. By complete reducibility of $V_d$
(being a finite-dimensional representation of the compact group
$K$) we can now use induction on the dimension of $V_d$.

\item{} This follows from the following commutative diagram (in which all the arrows are
isomorphisms of vector spaces and $k\in K$).

\[
\begin{CD}
{W}     @>>>  {W'}          @>{\pi (g)}>> {g W'}@={W'} @>>>   {W}\\
@V{d(k)}VV @V{\pi (k)}VV @V{\pi (g k g^{-1})}VV @V{\pi (k^g)}VV  @V{d(k^g)}V{{}^g d(k)}V\\
{W}     @>>>  {W'}          @>{\pi (g)}>> {g W'}@={W'} @>>>   {W}\\
\end{CD}
\]
\end{enumerate}

\end{proof}

\end{document}